\documentclass[11pt,letterpaper,reqno]{amsart}

\makeatletter%
\@mparswitchfalse%
\makeatother%

\reversemarginpar%

\usepackage{amsfonts}
\usepackage{amssymb}
\usepackage{amsmath}
\usepackage{amsthm}
\usepackage{mathtools}
\usepackage{amstext}
\usepackage{tikz-cd}

\usepackage[left=1.3in,right=1.3in,top=1.2in, bottom=1.3in]{geometry}

\definecolor{refkey}{HTML}{0003d6}
\definecolor{labelkey}{HTML}{083b19}

\usepackage[disable,colorinlistoftodos, textwidth=28mm, color=green!30]{todonotes}

\usepackage[ruled,noend,linesnumbered,algosection]{algorithm2e}
\usepackage{graphicx}

\definecolor{cit}{HTML}{117733}
\definecolor{lin}{HTML}{0003d6}
\definecolor{gree}{HTML}{6dc066}
\usepackage{hyperref}
\hypersetup{bookmarksnumbered,colorlinks=true,linkcolor=lin,citecolor=cit,urlcolor=black}

\usepackage[section]{placeins}
\usepackage{caption}

\newcommand{\mycomment}[1]{} %

\DeclareMathOperator{\dv}{div}
\DeclareMathOperator{\Pic}{Pic}
\DeclareMathOperator{\Cl}{Cl}
\DeclareMathOperator{\Div}{Div}
\DeclareMathOperator{\CDiv}{CDiv}

\DeclareMathOperator{\Span}{span}

\DeclareMathOperator{\Hom}{Hom}

\DeclareMathOperator{\PALP}{PALP}
\DeclareMathOperator{\GRD}{GRD}

\newcommand{\CC}{\mathbb{C}}
\newcommand{\Z}{\mathbb{Z}}

\newcommand{\R}{\mathbb{R}}
\newcommand{\PP}{\mathbb{P}}

\newcommand{\Q}{\mathbb{Q}}

\newcommand{\adj}{\leftrightarrow}
\newcommand{\nr}{\sim}
\DeclareMathOperator{\adjt}{adj}
\DeclareMathOperator{\near}{near}
\DeclareMathOperator{\nbd}{nbd}

\DeclareMathOperator{\AD}{AD}
\DeclareMathOperator{\conv}{conv}

\DeclareMathOperator{\relint}{relint}
\DeclareMathOperator{\Int}{int}

\newcommand{\vv}[1]{\lvert#1\rvert}
\newcommand{\vr}[1]{\langle#1\rangle}
\newcommand{\ra}{\to}

\newcommand{\mt}{\mapsto}

\newcommand{\txt}[1]{\text{ #1}} %
\theoremstyle{plain}%
\newtheorem{theorem}{Theorem}[section]

\newtheorem{corollary}[theorem]{Corollary}
\newtheorem{prop}[theorem]{Proposition}
\newtheorem{lemma}[theorem]{Lemma}

\theoremstyle{definition}%
\newtheorem{defi}[theorem]{Definition}

\newtheorem{example}[theorem]{Example}
\newtheorem{remark}[theorem]{Remark}

\newcommand{\pf}{{\em Proof}}%
\newcommand{\pfof}[1]{{\em Proof of #1}}%

\begin{document}
\title[Reflexive polytopes and Picard ranks of Gorenstein toric Fanos]{Reflexive polytopes and the Picard ranks of Gorenstein toric Fano varieties}
\author{Zhuang He}
\email{zhuang.he@unito.it}
\date{\today}

\begin{abstract}
We prove that the sum of the Picard ranks of a polar pair of Gorenstein toric Fano varieties of dimension $d\geq 3$ is at most the minimum of the number of facets and vertices of the corresponding pair of reflexive polytopes minus $(d-1)$. This is a generalization of Eikelberg's theory of affine dependences describing the Picard groups of toric varieties. The upper bound is achieved if and only if the polar pair is a simple-simplicial pair. 

\end{abstract}
\maketitle

\section{Introduction}\label{Intro}

Reflexive polytopes are lattice polytopes containing zero in the interior such that their polar polytope is also a lattice polytope. A lattice polytope is reflexive if and only if its polar is reflexive. Thus reflexive polytopes appear in pairs. Introduced by Batyrev \cite{Batyrev1994} for constructing toric Fano varieties with Calabi-Yau hypersurfaces that form mirror pairs, reflexive polytopes have been a central topic in mirror symmetry, convex geometry, and birational geometry. The toric variety associated with the spanning fan of a reflexive polytope is Gorenstein and Fano: the anti-canonical divisor class is Cartier and ample. The classes of reflexive polytopes up to lattice isomorphism one-to-one correspond to the isomorphism classes of Gorenstein toric Fano varieties.

A long-standing question on reflexive polytopes is to give upper bounds of the Picard ranks of the associated Gorenstein toric Fano varieties. There are only finitely many non-isomorphic reflexive polytopes in each dimension, while the numbers grow fast. For $\Q$-factorial toric Fano varieties of dimension $d$, the Picard ranks are at most $2d$ \cite{Casagrande2005}. A general upper bound for all Gorenstein toric Fano varieties remains unknown.

In the present paper we show a new approach to the question above. Let $P\subseteq N_\R$ be a reflexive polytope and $X_P$ be the toric variety from the spanning fan of $P$. We generalize Eikelberg's description \cite{Eikelberg1992} of the non-$\Q$-factorial toric Picard groups via spaces of affine dependences, and \cite[Cor. 2.18]{Eikelberg1993} that if $P$ is simple (that is, the polar of a simplicial polytope), then the Picard rank $\rho_X=1$. Our main result is:

\begin{theorem}\label{PQEintro}
Let $d\geq 3$. Let $X=X_P$ and $Y=Y_Q$ be two Gorenstein toric Fano varieties of dimension $d$, such that $P$ and $Q$ are a polar pair of reflexive polytopes. Then
\[\rho_X + \rho_Y \leq \min \{\vv{F(P)}, \vv{V(P)}\} -d +1,\]
where $\vv{F(P)}$ and $\vv{V(P)}$ are the number of facets and vertices of $P$ respectively, and $\rho_X, \rho_Y$ are the Picard ranks of $X$ and $Y$. 
Furthermore, the equality holds if and only if $P$ or $Q$ is simplicial, in which case $\rho_Y=1$ or $\rho_X=1$ respectively.
\end{theorem}

We remark that if two $d$-dimensional polytopes $P$ and $Q$ are polar to each other, then there is a one-to-one correspondence between the $k$-faces of $P$ and the $(d-1-k)$-faces of $Q$. In particular, $\vv{F(P)}=\vv{V(Q)}$ and $\vv{V(P)}=\vv{F(Q)}$. Thus the statement of Theorem \ref{PQEintro} is symmetric between $P$ and $Q$. 

To the best of the knowledge of the author, there have been no similar treatments in literature concerning the sum of the Picard ranks of a polar pair of non-$\Q$-factorial Gorenstein toric Fano varieties. In the simplicial-simple case, Theorem \ref{PQEintro} recovers \cite[Cor. 2.18]{Eikelberg1993}. To be precise, $P$ is simple if and only if $Q$ is simplicial, which is equivalent to $X_Q$ is $\Q$-factorial. Under these conditions, $\Pic(X_Q)_\Q=\Cl(X_Q)_\Q$, which has rank $\vv{V(Q)}-d$. On the other hand, by \cite[Cor. 2.18]{Eikelberg1993}, $\Pic(X_P)\cong \Z$. Our result particularly implies that $\vv{V(Q)}\leq \vv{F(Q)}$ for simplicial $Q$, which can be shown for example by the Lower Bound Theorem by \cite{Barnette1973}. In general, $\Pic(X_Q)_\Q$ is a subgroup of $\Cl(X_Q)_\Q$, which is of rank $\vv{V(Q)}-d$, while an upper bound of $\vv{V(Q)}$ is unknown for reflexive polytopes.

Assume $\vv{F(P)}\leq \vv{V(P)}$. That is, $P$ is the `simpler' one and $Q$ is the more `simplicial' one. We define 
\[\varepsilon = \varepsilon(P, Q) := \vv{F(P)}-d +1 - \rho_X- \rho_Y.\]
Then Theorem \ref{PQEintro} says $\varepsilon\geq 0$ for all reflexive polytopes of $d\geq 3$, and $\varepsilon= 0$ if and only if $(P,Q)$ is a simple-simplicial pair.

Even though Theorem \ref{PQEintro} is about a polar pair of reflexive polytopes, our proof is based on a combinatorial construction over a single reflexive polytope, which gives a new upper bound $\rho'=\rho'_P$ of the actual Picard rank $\rho_X$ of $X_P$. Given a reflexive polytope $P$, we consider all the permutations $L$ of the set of facets $F(P)$, and assign an integer $1\leq \rho_L\leq \infty$ to $L$: each facet $G$ in $L$ except the first one contributes to $\rho_L$ by one, if the affine span of those vertices of $G$ that are in some previous facets in $L$ is of dimension $d-2$, or $\infty$ if the dimension is smaller. Then we define the {\em facet complexity} $\rho'$ of $P$ as the minimum of all $\rho_L$ over all the permutations $L$ (Definition \ref{fc}).

Later we will refer to such permutations $L$ as {\em facet sequences}, or simply sequences over $F(P)$. It turns out that $\rho'$ is always a finite integer (Lemma \ref{finitecomplexity}). An important scenario is when $L$ is an {\em adjacent sequence} (Definition \ref{adjSeq}). That is, every subsequent facet is adjacent to some previous facet in $L$. It can be shown that every convex polytope allows an adjacent sequence of facets, hence $\rho'\leq \vv{F(P)}-1$.

\begin{theorem}[Theorem \ref{fc1}]\label{fc1intro}
Let $P$ be a reflexive polytope and $X=X_P$ be the associated Fano variety. Then
\[\rho_X \leq \rho' \leq\vv{F(P)}-1.\]

\end{theorem}

We describe the key steps proving $\rho_X \leq \rho'$. Eikelberg \cite{Eikelberg1992} showed that the group $\CDiv(X)_\Q$ of $\Q$-Cartier divisors on $X$ is isomorphic to the dual of a quotient vector space $H_\Q/\AD(\Sigma)$ (see Theorem \ref{ADCdiv}). The ambient space $H_\Q$ has a basis indexed by the vertices of $P$, thus $\dim H_\Q = \vv{V(P)}$. On the other hand, $\AD(\Sigma)$, the space of affine dependences, is a sum of subspaces $\AD(\sigma)$ over all the maximal dimensional cones $\sigma$ of the fan $\Sigma$ over the faces of $P$. The local structure of $\{\AD(\sigma)\}_\sigma$ enables us to construct a flag of subspaces of $H_\Q/\AD(\Sigma)$, from a sequence $L$ of facets. Each dimension increment in the flag is bounded by the contribution of the specific facet to $\rho_L$, and hence $\rho_L+d$ is an upper bound of $\dim H_\Q/\AD(\Sigma)$, which equals to $\rho_X+d$.

We remark that other generalizations of the space of affine dependences include \cite[Definition 3.5]{Fujita2019}.

Suppose we have a polar pair of reflexive polytopes $(P,Q)$. Facet sequences of $P$ and $Q$ determine $\rho'_P$ and $\rho'_Q$ respectively. We use the fact that a sequence of facets $T$ of $Q$ is dual to a sequence of vertices $T^*$ of $P$, and adjacency can be defined via duality to the latter.
With this, for any adjacent sequence $T$ consisting of all the facets of $Q$, we construct an adjacent sequence $L$ of facets of $P$, by inductively adding the facets containing the latest vertex added to $T^*$ in a compatible order. When $d\geq 3$, we estimate simultaneously the facet complexities $\rho'_X$ and $\rho'_Y$ by examining $\rho_{L}$ and $\rho_T$, which gives the estimation in Theorem \ref{PQEintro}.

The question on the Picard ranks of toric Fano varieties of dimension $d$ has been raised and studied since inception. Batyrev conjectured that smooth toric Fano varieties have Picard ranks at most $2d$ \cite{Ewald2012}. \cite{Voskresenskiĭ1985} first showed that the smooth toric Fano varieties have Picard ranks $\rho_X\leq 2d^2-d$. \cite[Thm. 8]{Debarre2000} improved this to an upper bound of growth rate $d^{3/2}$. \cite{Nill2005} proved that the upper bound $2d$ holds for simplicial reflexive polytopes whose polar polytope contains both a vertex $v$ and the point $-v$.
The conjecture $\rho_X\leq 2d$ is proved for all simplicial reflexive $P$ by \cite{Casagrande2005}. Equivalently, $\vv{V(P)}\leq 3d$ for a simplicial reflexive polytope. For the number of vertices, \cite[5.2]{Nill2005} conjectured that $\vv{V(P)}\leq 6^{d/2}$ for all reflexive $d$-polytopes. 

A related result on the upper bound of Picard ranks is the Generalized Mukai Conjecture \cite{Bonavero2003}, which is proved for simplicial toric Fano varieties by \cite{Casagrande2005}. \cite{Fujita2019} proved the conjecture for toric log Fano pairs and in particular, for all Gorenstein toric Fano varieties.

The paper is structured as follows. Notations on convex geometry, toric geometry, and preliminaries on reflexive polytopes are presented in Section \ref{Prelim}. Section \ref{PicardAFC} introduces the facet complexity and shows its relation with the actual Picard rank. We also discuss the algorithmic aspect of computing the facet complexity in this paragraph. Section \ref{PolarPair} proves Theorem \ref{PQEintro}.

{\bfseries Acknowledgment} We thank Cinzia Casagrande, Kento Fujita, Christian Haase, Johannes Hofscheier, Benjamin Nill, and Andrea Petracci for various communications and suggestions. The project relied on the following packages and software: {\em PALP} \cite{Kreuzer2004}, the {\em Graded ring database} \cite{GD,Kas08}, {\em SageMath} \cite{sagemath}, the packages \verb|NormalToricVarieties| \cite{NTVSource} and \verb|Polyhedra| \cite{PolyhedraSource, PolyhedraArticle} of {\em Macaulay2} \cite{M2}. The author has been funded by the European Union - Next Generation EU, Mission 4, Component 1, CUP D53D23005860006, within the project PRIN 2022L34E7W ``Moduli spaces and birational geometry''. 

\section{Preliminaries}\label{Prelim}
\subsection{Notations and facts in toric geometry and polyhedral geometry}
We refer to \cite{CLS2011} for a general reference of toric geometry. For conventions in convex geometry, we follow \cite{Nill2005} and \cite{Brøndsted1983}. Toric varieties in this paper are assumed to be over the complex numbers. Let $N:= \Z^d$ and $N_\R := N\otimes_\Z \R$. A toric variety $X=X_\Sigma$ is decided by a fan $\Sigma$ in $N_\R$. We assume $X$ is complete, which is equivalent to that the support $\vv{\Sigma}:=\bigcup_{\sigma\in \Sigma} \sigma$ of the fan $\Sigma$ is the whole vector space $N_\R$. In this case, $X$ is a $d$-dimensional variety.

Let $\Sigma(k)$ be the set of $k$-dimensional cones in $\Sigma$. 

We denote by $\conv(S)$ the convex hull of any set $S\subseteq N_\R$.
A lattice polytope $P$ in $N_\R$ is the convex hull of finitely many lattice points. The boundary of $P$ is denoted by $\partial P$, and the relative interior by $\relint(P)$, or simply $\Int(P)$ if $P$ is full-dimensional in $N_\R$. A $d$-dimensional polytope is often called a $d$-polytope.

A $k$-face of $P$ is a face of dimension $k$.
If $P$ is of dimension $m$, then a facet of $P$ is a $(m-1)$-dimensional face of $P$. An edge is a one-face of $P$, and is of the form $\overline{vw}:=\conv(v,w)$ for two vertices $v$ and $w$. Let $F(P)$ be the set of facets of $P$, and $V(P)$ the set of vertices of $P$.

The $f$-vector $(f_0,\cdots, f_m)$ is defined by letting $f_k$ be the number of $k$-faces of $P$. In particular, $f_0=\vv{V(P)}$ and $f_{m}=\vv{F(P)}$.

Let $F$ be a facet of $P$. The normal ray of $F$ has rational slopes, hence contains a lattice point. We denote by $u_F$ the shortest nonzero inner normal vector of $F$ that is a lattice point.

In the following we always assume $P$ is full dimensional in $N_\R$ and $0\in \Int(P)$. The spanning fan $\Sigma(P)$ is the fan in $N_\R$ consisting of the cones over the faces of $P$. We write $X_P$ for the toric variety associated with the spanning fan $\Sigma(P)$. Then $X_P$ is projective. Conversely, any projective toric variety is constructed out of the spanning fan of a polytope $P\in N_\R$ such that $0\in \Int(P)$ (see \cite[Example 3]{Debarre2000}).

The toric variety $X_P$ is smooth if and only if $P$ is a smooth polytope, i.e., for any facet $F$ of $P$, the vertices of $F$ is part of a basis for $N$. Similarly, $X_P$ is $\Q$-factorial if and only if $P$ is simplicial, i.e., any facet $F$ of $P$ is a simplex. 

Back to the polyhedral geometry of $P$. We denote by $v,w,x,y,z,t,\cdots$ the vertices, and $F, G, H, K,\cdots$ the facets of $P$. The valence $e(v)$ of a vertex $v$ is the number of distinct facets containing $v$. 

The polar, or dual polytope of $P$ in $N_\R$ is
\[ P^* := \{x\in M_\R: \vr{x,y}\geq -1 \txt{for all } y\in P\}.\]
where $M:=\Hom(N,\Z)$ is the dual lattice of $N$. There is a one-to-one correspondence between the $k$-faces of $P$ and $(d-1-k)$-faces of $Q$. A polytope is called {\em simple} if $P^*$ is simplicial. We note the following facts:
\begin{lemma}\label{simple}
\cite[Thm. 12.11, 12.12]{Brøndsted1983}
Let $P$ be a $d$-polytope. The following are equivalent:
\begin{enumerate}
\item $P$ is simple;
\item all the vertices of $P$ are of valence $d$;
\item all the vertices of $P$ are incident to exactly $d$ edges of $P$.
\end{enumerate}
\end{lemma}

We consider the relative positions between facets and vertices of $P$. We say two distinct facets $F$ and $G$ are {\em adjacent} if $F\cap G$ is a $(d-2)$-face of $P$. We say $F$ and $G$ are {\em near} if $F\cap G \neq\emptyset$. Taking the dual, we say two distinct vertices $v$ and $w$ are adjacent if the line segment $\overline{v w}$ is an edge of $P$, and near if $v$ and $w$ are in a common facet. Clearly, adjacent facets or vertices are near. We write $F\adj G$, and $v\adj w$ for adjacency, and $F\nr G$, and $v\nr w$ for nearness. 

We also write $F\nr^v G$ if $v\in F\cap G$, and $v\nr^F w$ if $v\in F$ and $w\in F$. Finally, any $(n-2)$-dimensional face $\tau$ of $P$ is contained in exactly two facets $F$ and $G$, and $\tau= F\cap G$. In this case we write $F\adj^\tau G$.

For a facet $F$ or vertex $v$, we let $\adjt(F), \near(F)$ or $\adjt(v), \near(v)$ be the set of facets or vertices that are adjacent or near $F$ or $v$. For a vertex $v$, define $\nbd(v)$ as the set of facets containing $v$.

Given a convex $d$-polytope $P$ with $d\geq 2$, the vertices $V(P)$ and edges $E(P)$ naturally form an undirected graph $\mathcal{G}(P)$.

\begin{lemma}\label{graphlemma}\cite[Thm. 15.5]{Brøndsted1983}
Suppose $F$ is a facet of $P$, and $M$ is a possibly empty subset of $V(F)$. Then the subgraph of $\mathcal{G}(P)$ spanned by $V(P)\backslash M$ is connected.
In particular, $\mathcal{G}(P)$ is connected.
\end{lemma}

\subsection{Reflexive polytopes}
A lattice polytope $P$ is called reflexive if $0\in \relint(P)$ and $P^*$ is also a lattice polytope. Equivalently, $P$ is reflexive if and only if $\vr{u_F, x} = -1$ for any facet $F$ and any boundary point $x\in F$.

Note $P$ is reflexive if and only if $P^*$ is reflexive. If $P$ is reflexive, then $0$ is the unique lattice point in $\relint(P)$.

A normal variety $X$ over $\CC$ is {\em Gorenstein} if $K_X$ is Cartier. $X$ is Gorenstein and Fano if $-K_X$ is an ample Cartier divisor. The toric variety $X_P$ is Gorenstein and Fano if and only if $P$ is reflexive \cite[Thm. 4.19]{Batyrev1994}. The isomorphism classes of reflexive polytopes 1-1 correspond to the isomorphism classes of Gorenstein toric Fano varieties.

There are only finitely many non-isomorphic reflexive polytopes in each dimension $d$. For $d=1$ the unique reflexive polytope is the line interval $[-1,1]$. There are $16$ reflexive polygons (see \cite[Prop. 4.1]{Nill2005M} for a list). For $d=3$ and $4$, there are $4,319$ and $473,800,776$ reflexive polytopes respectively, by a monumental work by Kreuzer and Skarke \cite{Kreuzer1997,Kreuzer1998,Kreuzer2000a,Kreuzer2000b}. Table (\ref{Tenum}) lists all the counting relevant in this paper, where the number of smooth $4$-polytopes is by \cite{Batyrev1999} and \cite{Sato2000}, and the others by Kreuzer and Skarke.

The complete classification by Kreuzer and Skarke is stored in the {\em PALP} database \cite{Kreuzer2004}, and can be accessed via packages in {\em SageMath} \cite{sagemath} and {\em Macaulay2} \cite{M2}. For the $4,319$ reflexive $3$-polytopes, we define the PALP ids of them as their positions in the list of PALP database, ranging from $0$ to $4,318$. The same polytopes are also indexed in the Graded Ring Database \cite{GD,Kas08} with different ids, which we refer to as the GRD ids, from $1$ to $4,319$. We write $\PALP(n)$ or $\GRD(n)$ for the polytope with PALP id $n$, or GRD id $n$. For example, $\PALP(0)=\GRD(1)$, and is the smooth simplex corresponding to $\PP^3$; $\PALP(88)=\GRD(442)$, with $5$ vertices, and the representative in the PALP database has vertices $(1,0,0), (0,1,0), (0,0,1), (-4,-2,-1)$ and $(-2,-2,1)$.

\begin{table}
\centering
\caption{Counting reflexive polytopes in low dimensions}\label{Tenum}
\begin{tabular}{|c | r r r|} 
 \hline
 Dimensions & reflexive & simplicial & smooth \\ [0.5ex] 
 \hline
 2 & 16 & 16 & 5 \\
 3 & 4,319 & 194 & 18 \\
 4 & 473,800,776 & & 124 \\
 \hline
\end{tabular}
\end{table}

\section{Picard groups, affine dependences, and the facet complexity}\label{PicardAFC}
\subsection{Affine dependences}

It is a standard result in toric geometry, that for a toric variety $X=X_\Sigma$, the group $\Div_T(X)$ of the $T$-invariant Weil divisors over $X$ is freely generated by the classes of the torus-invariant divisors $D_\sigma$, where $\sigma$ runs over the $1$-dimensional faces $\Sigma(1)$. When $X$ is complete, the map $M\ra \Div_T(X)$ sending each $m$ to the principle divisor $\dv(\chi^m)$ is injective, so $\Cl(X)_\Q\cong \Q^{\vv{\Sigma(1)}-d}$. When $-K_X$ is $\Q$-Cartier, there is a lattice polytope $P$ in $N_\R$ whose vertices are the ray generators of $\Sigma$ \cite[\S 2.4]{Debarre2000}. Therefore, the rank of the $\Q$-class group equals to $\vv{V(P)}-d$. If $X$ is in addition $\Q$-factorial, or equivalently, if $P$ is simplicial, then the Picard rank $\rho_X$ of $X$ equals to $\vv{V(P)}-d$ too. 

In general $\Cl(X)$ may have torsions, while $\Pic(X_P)$ is a free abelian group of rank $\rho_X\leq \vv{V(P)}-d$ (see \cite[Prop. 4.2.5]{CLS2011}).

\begin{theorem}\cite{Casagrande2005}
If $P$ is a simplicial reflexive $d$-polytope. Then $\vv{V(P)}\leq 3d$, and $\rho_X\leq  2d$. The equality holds if and only if $d=2k$ is even, and $X$ is the product $(S_3)^{k}$, where $S_3$ is the del Pezzo surface of degree $6$.
\end{theorem}

We recall the description of the Picard group for an arbitrary toric variety by Eikelberg \cite{Eikelberg1992}. Assume $X=X_\Sigma$ is a complete normal toric variety. Let $m:=\vv{\Sigma(1)}$, and $D_i$ be the torus-invariant divisor corresponding to the ray $\tau_i\in \Sigma(1)$, for $i=1,\cdots, m$. Let $v_i$ be the ray generator of $\tau_i$. Fix a $T$-invariant $\Q$-Cartier divisor $L=l_1 D_1 + \cdots + l_m D_m$ such that all $l_i\neq 0$. Define $H:= \Z^{\oplus \Sigma(1)}\cong \Z^m$, with the standard basis $h_1,\cdots, h_m$.

\begin{defi}\label{ad}\cite[Thm. 3.2]{Eikelberg1992}
For a fixed $L$ defined as above, the space of affine dependences of a maximal dimensional fan $\sigma\in \Sigma(d)$ is the linear subspace
\begin{equation}
\AD(\sigma):=\left\{\sum_{i=1}^m a_i h_i\in H_\Q: 
\begin{array}{l}
a_i =0 \txt{if } \tau_i \not\subseteq \sigma,\\ \sum_{i=1}^m a_i \left(\displaystyle\frac{v_i}{l_i}\right)= 0, \txt{and } \sum_{i=1}^m a_i =0.
\end{array}
\right\} \subseteq H_\Q.
\end{equation}

Define the space of {\em affine dependences} of $X$ (or $\Sigma$) as 
\begin{equation}
\AD(X_\Sigma)=\AD(\Sigma):= \sum_{\sigma\in \Sigma(d)} \AD(\sigma) \subseteq H_\Q.
\end{equation}
\end{defi}

\begin{lemma}\cite{Eikelberg1992}\label{ECartierCenter}
Assume $X$ is a complete toric variety. Then there exists a $T$-invariant $\Q$-Cartier divisor $L=l_1 D_1 + \cdots + l_m D_m$ such that all $l_i\neq 0$. 
\end{lemma}

Let $\CDiv_T(X)_\Q$ be the group of $T$-invariant Cartier divisors on $X$.
\begin{prop}[\cite{Eikelberg1992}, also see {\cite[Thm. 3.1]{Fujita2019}}]\label{ADCdiv}
There is a perfect pairing:
\begin{align}
\CDiv_T(X)_\Q \times H_\Q/\AD(\Sigma)&\ra \Q\\
\left(D=\sum_{i=1}^m r_i D_i, \;\quad\sum_{i=1}^m a_i h_i\right) &\mt \sum_{i=1}^m \left( \frac{r_i}{l_i}\right) a_i.
\end{align}
\end{prop}
\pf. Consider the pairing $\vr{\cdot,\cdot}:\CDiv_T(X)_\Q \times H_\Q\ra \Q$ defined by the mapping above. The pairing $\vr{\cdot,\cdot}$ is indeed well-defined. Since $\vr{D,h_i}=r_i/l_i$, if $\vr{D,h_i}=0$, then $r_i=0$. Therefore if $\vr{D,h_i}=0$ for all $h_i$, then $D=0$.

By \cite[Claim A, page 516]{Eikelberg1992}, if $\sum_{i=1}^m (r_i/l_i)\cdot a_i = 0$ for all $\sum_{i=1}^m a_i h_i\in \AD(\Sigma)$, then $\sum_{i=1}^m r_i D_i$ is $\Q$-Cartier. On the other hand, by \cite[Claim B, page 516]{Eikelberg1992}, if $\sum_{i=1}^m r_i D_i$ is $\Q$-Cartier, then $\sum_{i=1}^m (r_i/l_i)\cdot a_i = 0$ for all $\sum_{i=1}^m a_i h_i\in \AD(\Sigma)$. Therefore the pairing factors through $\AD(\Sigma)$, and the induced pairing is perfect. \qed

\begin{corollary}\cite[Theorem 3.2]{Eikelberg1992}\label{EikelbergT}
The Picard group $\Pic(X)$ is a free abelian group of rank $m-d-\dim \AD(\Sigma)$.
\end{corollary}

\begin{remark}
We have $\AD(\sigma)=0$ if and only if $\sigma\in \Sigma(d)$ is a simplex. As a result, $\AD(\Sigma)=0$ if and only if $\Sigma$ is simplicial.
\end{remark}

For each $i$, let $v'_i:= v_i/l_i$. We call $\{v'_i\}$ the scaled ray generators. In general $v'_i\not\in\Z$.

\begin{remark}
The space $H_\Q$ is canonically isomorphic to the space $\tilde{A}_1(X)_\Q$ spanned by the $T$-invariant curves of $X$, via $h_i \mt V(\tau_i)$, where $\tau_i\in \Sigma(1)$. The perfect pairing between $\Div_T(X)$ and $\tilde{A}_1(X)_\Q$ is compatible with the pairing in Theorem \ref{ADCdiv}. Under this view, $\AD(\sigma)$ can be identified with linear relations among the scaled ray generators $\{v'_i\}$ that vanish on all the $T$-invariant Cartier divisors of $X$.
\end{remark}

We examine the subspaces of $H_\Q/\AD(\Sigma)$ spanned by the scaled ray generators in a maximal dimensional cone. For any generator $h_i$ of $H$, we write $\overline{h}_i$ for the class of $h_i$ in $H_\Q/\AD(\Sigma)$. Recall that $t$ points $x_1,\cdots, x_t$ in a vector space $V$ over $k$ are called {\em affinely dependent} if there exist numbers $a_i\in k, i=1,\cdots, t$, such that $\{a_i\}$ are not all zero, $\sum_i a_i=0$, and $\sum_i a_i t_i =0$. The affine subspace spanned by $\{x_i\}_{i=1}^t$ is $V=\{\sum_{i=1}^t a_i x_i: \sum_{i=1}^t a_i=0\}$, and the affine dimension of $\{x_i\}_{i=1}^t$ is the dimension of $V$ as a vector space.

\begin{lemma}\label{adad}
Suppose $\tau\in \Sigma(k)$, and $\{v_{i_1},\cdots,v_{i_s}\}$ is a subset of ray generators of $\tau$. Then $v'_{i_1},\cdots,v'_{i_s}$ are affinely independent in $N_\Q$ if only only if $\overline{h}_{i_1},\cdots, \overline{h}_{i_s}$ are linearly independent in $H_\Q/\AD(\Sigma)$.
\end{lemma}
\pf. Without lost of generality, we can relabel $\{v_i\}$ and assume $i_j=j$ for $j=1,\cdots, s$.

``$\Longrightarrow$": Suppose $\overline{h}_1,\cdots, \overline{h}_s$ are linearly dependent. Then there exist $a_i\in \Q, i=1,\cdots, s$, not all zero, such that $\sum_{i=1}^s a_i h_i \in \AD(\Sigma)$. Note that for every $\sigma\in \Sigma(d)$, if $\sum_{i=1}^m b_i h_i \in \AD(\sigma)$, then $\sum_{i=1}^m b_i (v_i/l_i)= 0$, and $\sum_{i=1}^m b_i =0$. Therefore by linearity, for every $\sum_{i=1}^m b_i h_i \in \AD(\Sigma)$, $\sum_{i=1}^m b_i (v_i/l_i)= 0$, and $\sum_{i=1}^m b_i =0$ too. In particular, if we choose $b_i=a_i$ for $i=1,\cdots, s$, and $b_i=0$ for $i>s$, then $\sum_{i=1}^m b_i h_i =\sum_{i=1}^s a_i h_i\in \AD(\Sigma)$. Thus
\[\sum_{i=1}^s a_i \cdot \frac{v_i}{l_i}= 0, \txt{and } \sum_{i=1}^s a_i =0.\]
Therefore $\{v'_1,\cdots,v'_s\}$ are affinely dependent.

``$\Longleftarrow$": Suppose $\{v'_1,\cdots,v'_s\}$ are affinely dependent. Then there are $a_i\in \Q, i=1,\cdots, s$, not all zero, such that $\sum_{i=1}^s a_i (v_i/l_i)= 0$, and $\sum_{i=1}^s a_i =0$. There exists $\sigma\in \Sigma(d)$ containing $\tau$, so by the assumption, $v_i\in \sigma$ for all $i$. Let $b_i=a_i$ for $i=1,\cdots, s$, and $b_i=0$ for $i>s$. By Definition \ref{ad}, we have $\sum_{i=1}^s b_i h_i =\sum_{i=1}^m b_i h_i \in \AD(\sigma)\subseteq \AD(\Sigma)$. That is, $\overline{h}_1,\cdots, \overline{h}_s$ are linearly dependent in $H_\Q/\AD(\Sigma)$.\qed

\begin{remark}
The proof above shows that the ``only if'' part of the lemma holds even when $\{v_i\}$ are not in a common cone. 
\end{remark}

\begin{lemma}\label{fixedDim}
Suppose $\tau\in \Sigma(k)$, and $v_{i_1},\cdots,v_{i_s}$ are all the distinct ray generators of $\tau$. Suppose $v'_{i_1},\cdots,v'_{i_s}$ span an affine subspace of dimension $k'$. Then $\overline{h}_{i_1},\cdots, \overline{h}_{i_s}$ span a linear subspace of dimension $k'+1$. 

\end{lemma}
\pf. By Lemma \ref{adad}, $\{v'_{i_j}\}_j$ are affinely independent if and only if $\{\overline{h}_{i_j}\}_j$ are linearly independent in $H_\Q/\AD(\Sigma)$, so the two dimensions have a difference of $1$. \qed

In later paragraphs, we will only consider the case where $L=-K$ and is Cartier. In other words, $X$ is a Gorenstein toric Fano variety, and $X=X_P$, where $P$ is the reflexive polytope in $N_\R$ with vertices $\{v_i: i=1, \cdots, m\}$. Since $-K\sim \sum_{i=1}^m D_i$, we find that under such assumptions, all $l_i=1$, so the assumption that all $l_i\neq 0$ is satisfied. We can restate the lemmas \ref{adad} and \ref{fixedDim} in terms of facets and vertices of $P$. In particular, Lemma \ref{fixedDim} says that if the vertices $v_{i_1},\cdots,v_{i_s}$  are in a common facet, and span an affine subspace of dimension $k'$, then $\overline{h}_{i_1},\cdots, \overline{h}_{i_s}$ span a linear subspace of dimension $k'+1$.

\subsection{Facet complexity and adjacent sequences}
Fix a reflexive polytope $P\subseteq N_\R$ of dimension $d$. Let $m:=\vv{F(P)}$. We frequently consider permutations of a set $\mathcal{F}$ of facets of $P$, that is, a sequence consisting of all the facets in $\mathcal{F}$ without repetition. We call such sequence a facet sequence of $P$ over $\mathcal{F}$, or simply a sequence over $\mathcal{F}$. Similarly we define (vertex) sequences over a set of vertices $\mathcal{V}$.

\begin{defi}\label{adjSeq}
\begin{enumerate}
\item 
We say a facet sequence $L=(G_0,\cdots, G_{s-1})$ of length $s$ is an {\em adjacent sequence}, if either $s=1$, or for each $i, i=1\cdots, s-1$, there exists $r_i\leq i-1$ such that $G_{r_i}\adj G_i$.
\item 
We say a vertex sequence $T=(v_0,\cdots, v_{s-1})$ of $P$ is an {\em adjacent sequence}, if the dual sequence $T^*$ is an adjacent facet sequence of $P^*$. Equivalently, $T$ is an adjacent sequence if either $s=1$, or for each for each $i, i=1\cdots, s-1$, there exists $r_i\leq i-1$ such that $v_{r_i}\adj v_i$.
\end{enumerate}
\end{defi}

To any sequence $L:=(G_0, G_1,\cdots, G_{m-1})$ over $V(P)$, we can assign a value $\rho_L\in \Z_{\geq 0} \cup \{\infty\}$ as follows: for each $i\geq 0$, let $S_i:= \cup_{j=0}^i V(G_j)$ be the set of vertices in any of $G_0,\cdots, G_i$. Let $w_i$ be the affine dimension of the points in $S_{i-1}\cap V(G_i)$ if $S_{i-1}\cap V(G_i)\neq\emptyset$, and set $w_i=-1$ otherwise. For each $i\geq 1$, define $a_i:=d-1-w_i$ if $w_i\geq d-2$, and $\infty$ otherwise. Since the vertices in $G_i$ are on the same affine hyperplane $\Span(G_i)$, we must have $w_i\leq d-1$. In other words, the definition of $a_i$ is equivalent to:
\begin{align*}
a_i := 
\begin{cases}
	0, & \txt{if } w_i = d-1;\\
	1, & \txt{if } w_i = d-2;\\
	\infty, & \txt{else}.
\end{cases}
\end{align*}
Lastly, let $\rho_L:=\sum_{i=1}^{m-1} a_i$, where the convention is $\infty + n = \infty$ for $n\in \Z \cup \{\infty\}$.

\begin{defi}\label{fc}
	Let $P$ be a reflexive polytope of dimension $d$. We define the {\em facet complexity} $\rho'=\rho'(P)$ of $P$ to be the minimum of $\rho_L$ for all the facet sequences $L$ over $F(P)$.
\end{defi}

Before giving an example, we first show that the complexity $\rho'$ is always a finite number. 

\begin{lemma}\label{adjSeqE}
Let $P$ be a convex polytope. Then the following hold:
\begin{enumerate}
\item There exists an adjacent sequence over $F(P)$.
\item Fix a vertex $v$ of $P$. If $\mathcal{F}$ is a subset of all the facets that contain $v$, then there exists an adjacent sequence $L$ over $F(P)\backslash \mathcal{F}$.

\item Any sequence over a subset $\mathcal{F}=(F_1,\cdots, F_t)$ of $F(P)$ can be extended to a sequence $L'$ over $F(P)$, such that for each $i\geq t+1$, $F_i\adj F_j$ for some $j<i$.
\item In particular, any adjacent sequence $L$ over a subset $\mathcal{F}\subseteq F(P)$ can be extended into an adjacent sequence $L'$ over $F(P)$.
\end{enumerate}
\end{lemma}
\pf. By duality we need only prove the dual statements for vertices. 
First we prove (3). Let $\mathcal{F}^*$ be any vertex sequence of $Q=P^*$. If $\mathcal{F}^*=\emptyset$, then we can choose any vertex $v$ as the initial vertex of $(L')^*$. Otherwise, suppose $\mathcal{F}\subsetneq F(P)$. By Lemma \ref{graphlemma}, the graph $\mathcal{G}(Q)$ is connected, so there exists a vertex $v_0\not\in\mathcal{F}^*$ such that $v_0\adj v'$ for some $v'\in \mathcal{F}^*$. Inductively, we can keep adding vertices with an edge to those previous vertices, until no such external edges exist. When this procedure terminates, all the vertices must have been added to the sequence $(L')^*$, whose dual $L'$ is the targeted sequence. This proves (3). For (2), the dual $\mathcal{V}$ of $\mathcal{F}$ is a subset of $V(F_v)$ on $Q$, where $F_v$ is the dual facet of $v$. Consider $\mathcal{V}':=V(P)\backslash \mathcal{V}$. By Lemma \ref{graphlemma}, the subgraph of $\mathcal{G}(Q)$ over $\mathcal{V}'$ is connected. Hence it holds by the same proof of (3). (4) is a corollary of (3), noting that if $L$ is an adjacent sequence, then so is $L'$. Finally, (1) is a special case of (4) when $L$ is empty.\qed

\begin{lemma}\label{finitecomplexity}
If $P$ is a reflexive polytope, then $\rho'<\infty$.
\end{lemma}
\pf. By Lemma \ref{adjSeqE}(1), there exists an adjacent sequence $L=(G_0, G_1,\cdots, G_{m-1})$ over $F(P)$. For $L$, each $a_i<\infty$, since $S_{i-1}\cap V(G_i)$ contains all the vertices on a facet of $G_i$, so their affine span has affine dimension at least $d-2$. Thus $\rho_L<\infty$, and $\rho'<\infty$.\qed

\begin{example}\label{4309}
Consider the $3$-dimensional reflexive polytope $P=\PALP(4309)=\linebreak\GRD(2356)$. This $P$ is the unique $3$-dimensional reflexive polytope with precisely $14$ vertices and $12$ facets. We label the facets from $0$ to $11$ as in Figure \ref{4309figure}. The adjacent sequence $L_1=(0,1^*,7^*,11^*,5^*,9,3,8,6^*,2,4,10)$ has $\rho_{L_1}=5$, where facets marked with ($*$) are those with $a_i=1$. However, there are other sequences of lower $\rho_L$. For example, the adjacent sequence $L_2=(0,1^*,2,3^*,4,5,6,7,8,9,10,11)$ has $\rho_{L_2}=2$. 

In fact it is easy to see that $\rho'= 2$. In dimension $3$, let $k_i:=\vv{S_{i-1}\cap V(G_i)}$. Then $a_i=0$ if and only if $k_i\geq 3$, and $a_i=1$ if and only if $k_i=2$. Suppose a sequence $L$ has $\rho_L=1$. Then $L$ must satisfy $G_1\adj G_0$, and $a_1=1$. Then it is impossible to arrange all the remaining facets so that all the $a_i=0$ for $i\geq 2$. This can also be verified with a computer program (see Section \ref{algorithms}). 

\begin{figure}[h]

\includegraphics[width=0.8\textwidth]{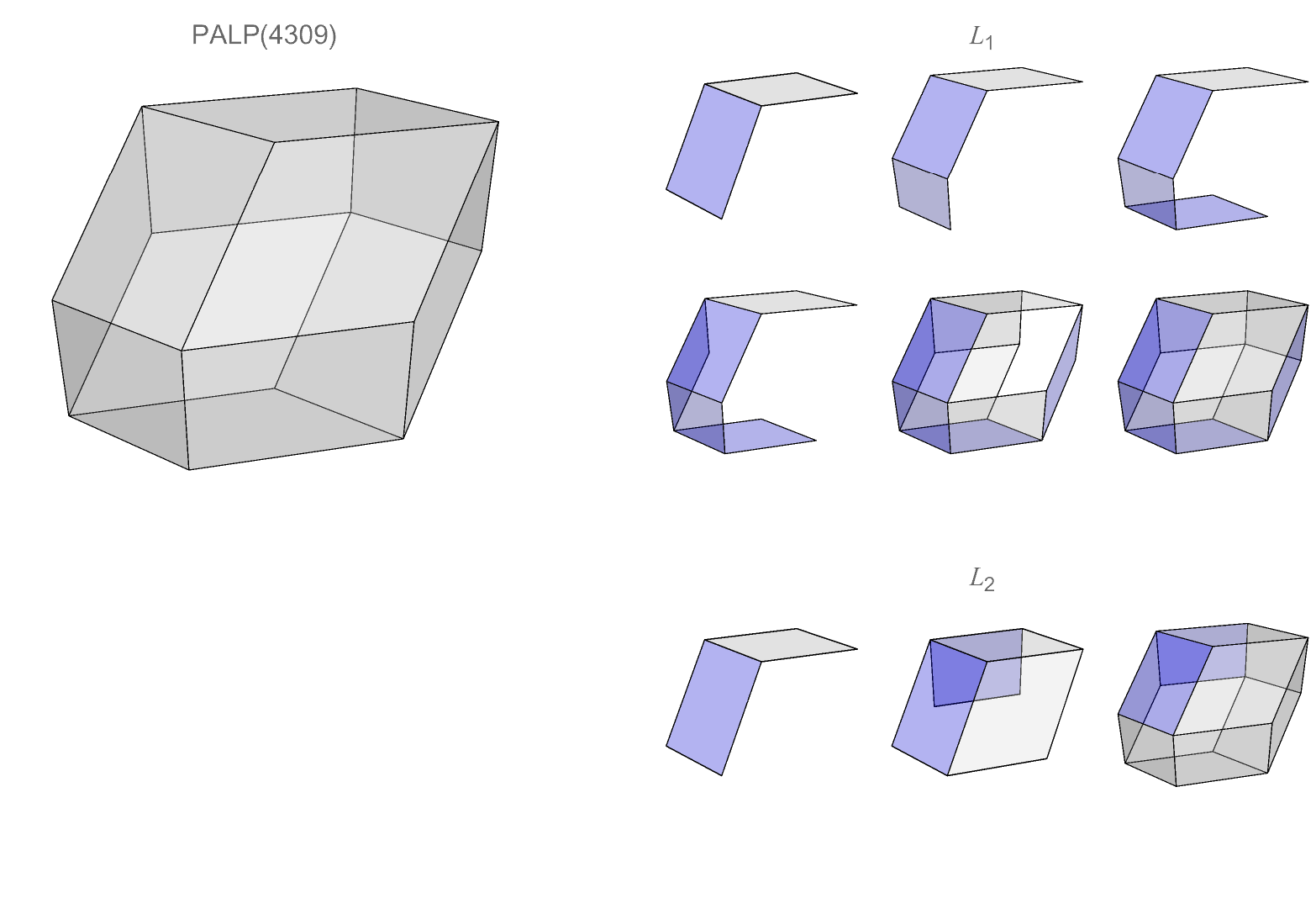}
\caption{\small The facet complexity of the reflexive polytope of PALP id 4309. Left: facet labels. Top: $0$; front: $1,2,7,8$ (left to right, then top to bottom); left: $5$; right: $6$; back: $3,4,9,10$; bottom: $11$. Top right: the sequence $L_1$ with $\rho_{L_1}=5$; bottom right: $L_2$ with $\rho_{L_2}=2$. Blue facets are those with (*) and $a_i=1$.}\label{4309figure}
\end{figure}
\end{example}

We can compute the Picard rank for the example \ref{4309} above, and it turns out $\rho_X=2$ too. Our main result in this paragraph shows the relations between $\rho_X$ and $\rho'$. 

\begin{remark}
If $Q$ is the polar of $P$, then the Picard rank of a projective toric variety $X_Q$ also equals to the dimension of the cone of Minkowski summands of the polytope $P$ (\cite[Thm. 4.16]{Eikelberg1993}, also see \cite[\S 6]{CLS2011}). For example, the polytope in Example \ref{4309} is a zonotope, and it equals to the Minkowski sum of $4$ line segments. \footnote{We thank Christian Haase for pointing this out.} Then one can show that the cone of Minkowski summands of $P$ is of dimension $4$, so $\rho_X = 4$.
\end{remark}

\begin{theorem}\label{fc1}
	Let $P$ be a reflexive polytope in $N_\R$. Let $X$ be the associated toric Fano variety of $P$. Let $\rho$ be the rank of $\Pic(X)$ and $\rho'$ be the facet complexity of $P$. Then
	\[\rho \leq \rho'\leq \vv{F(P)}-1.\]
\end{theorem}

\pf. By Lemma \ref{finitecomplexity}, $\rho'<\infty$. In particular, by Lemma \ref{adjSeqE}, there exists an adjacent sequence $L$ over $F(P)$. Then $\rho_L<\infty$, $\rho_L\leq \vv{F(P)}-1$, since each $a_i\leq 1$. Hence $\rho'\leq \vv{F(P)}-1$.

So we need only show that $\rho\leq \rho'$. For any subset $S=\{v_{i_1},\cdots, v_{i_t}\}$ of $V(P)$, let $H(S)$ be the subspace spanned by the corresponding generators $\{\overline{h}_{i_j}\}: 1\leq j\leq t$ in $H_\Q/\AD(\Sigma)$. Now we need only show that for any sequence $L=(G_0, G_1,\cdots, G_{m-1})$, $\rho\leq \rho_L$. For this $L$ and for each $i\geq 1$, by the dimension formula we have:
\begin{align*}
\dim H(S_{i-1}\cup V(G_i)) &\leq \dim H(S_{i-1}) + \dim H(V(G_i)) - \dim H(S_{i-1})\cap H(V(G_i))\\
									&\leq \dim H(S_{i-1}) + \dim H(V(G_i)) - \dim H(S_{i-1}\cap V(G_i)),
\end{align*}
where we use that $H(S_{i-1}\cap V(G_i))\subseteq H(S_{i-1})\cap H(V(G_i))$ by definition. 
\mycomment{
Suppose $S_1, S_2$ are two subsets in $V(P)$. 
\begin{align*}
\dim H(S_1\cup S_2) &= \dim H(S_1) + \dim H(S_2) - \dim H(S_1)\cap H(S_2)\\
									&\leq \dim H(S_1) + \dim H(S_2) - \dim H(S_1\cap S_2),
\end{align*}
where we use that $H(S_1\cap S_2)\subseteq H(S_1)\cap H(S_2)$ by definition. Now consider any sequence $L=(G_0, G_1,\cdots, G_{m-1})$. For each $i\geq 1$, Apply the formula above to $S_{i-1}$ and $V(G_i)$, we find:
\[\dim H(S_{i-1}\cup V(G_i)) \leq \dim H(S_{i-1}) + \dim H(V(G_i)) - \dim H(S_{i-1}\cap V(G_i)).\]
}
Note that $S_i= S_{i-1}\cup V(G_i)$. Since $S_{i-1}\cap V(G_i)\subseteq G_i$, we can apply Lemma \ref{fixedDim} to $V(G_i)$ and $S_{i-1}\cap V(G_i)$. Since $V(G_i)$ span the affine hyperplane containing $G_i$, $\dim H(V(G_i))= d$. The vertices $S_{i-1}\cap V(G_i)$ have affine dimension $w_i$ if nonempty, so $\dim H(S_{i-1}\cap V(G_i))=w_{i}+1$. Otherwise, $H(S_{i-1}\cap V(G_i))=\emptyset$ and $w_i=-1$, so the equality also holds. Therefore $\dim H(S_i) \leq \dim H(S_{i-1}) + d - (w_i+1)$. That is, 
\begin{equation}\label{eqHS}
\dim H(S_i) - \dim H(S_{i-1}) \leq a_i,
\end{equation}
since $a_i \geq d - (w_i+1)=d-1-w_i$. By Lemma \ref{fixedDim} again, $\dim H(S_{0}) = d$. Note $H(S_{m-1}) = H_\Q/\AD(\Sigma) $. Summing up (\ref{eqHS}) for $i=1,\cdots, m-1$, we have
\[\rho =\dim H_\Q/\AD(\Sigma) -d\leq \sum_{i=1}^{m-1} a_i =\rho_L.\]
\qed

\begin{remark}
The proof above of Theorem \ref{fc1} shares essentially the same idea with the proof of \cite[Theorem 2.16]{Eikelberg1993}.
\end{remark}

\begin{remark}
Our main theorem \ref{PQEintro} is an improved estimation on $\rho'$, and the bound $\vv{F(P)}-1$ here is never sharp for $d\geq 3$. In fact, when $P$ is a simple polytope of dimension $d\geq 3$, Eikelberg \cite{Eikelberg1993} proved that $\rho_{X} = 1$.
\end{remark}

\begin{example}\label{exsim}
	If $P$ is simplicial, then $\rho=\rho'= \vv{V(P)} -d$. Firstly, $X$ is $\Q$-factorial, so $\Pic(X)_\Q =\Div(X)_\Q$, and $\rho=\vv{V(P)} -d$. Since $P$ is simplicial, every facet of $P$ is a $(d-1)$-simplex. For any sequence $L=(G_0, G_1,\cdots, G_{m-1})$ over $F(P)$, $a_i=1$ if and only $V(G_i)\backslash S_{i-1}$ contains exactly one vertex $v_i$, and $a_i=0$ if and only $V(G_i)\subseteq S_{i-1}$. Therefore, if $\rho_L<\infty$, then $\rho_L=\sum_{i=1}^{m-1} a_i=\vv{V(P)} - d$.
\end{example}

\begin{example}\label{neqpr}
	In dimension two, all the $16$ reflexive polytopes are simplicial, so $\rho'=\rho$, and equals to the number of their edges minus $2$. When $d=3$, there are examples where $\rho\neq \rho'$. We implemented the Algorithm \ref{recursive alg} in {\em Macaulay2} and run it in dimension $3$. Among the $4,319$ reflexive polytopes of dimension $3$, $\rho= \rho'$ holds for all except $23$ polytopes. These $23$ polytopes consist of $22$ with $\rho'=2$ and $\rho=1$, and a unique one with $\rho'=3$ and $\rho=2$. The PALP ids of the $23$ exceptional ones are: 3385, 3386, 3522, 3523, 3847, 3848, 4000, 4015, 4016, 4045, 4046, 4071, 4072, 4083, 4084, 4096, 4097, 4252, 4253, 4281, 4316, 4317.
	The only one with $\rho'=3$ is $\PALP(4252)$.
	
It remains unknown whether there are reflexive polytopes with $\rho< \rho'$ in higher dimensions.
\end{example}

\subsection{Computing the facet complexity}\label{algorithms} 

Here we present an algorithm computing the facet complexity. Results here are not needed in later proofs, and reader can skip this paragraph for the main results.

From a computational point of view, the facet complexity $\rho'$ can be computed via a depth-first search, while it is better interpreted as a problem of dynamic programming, and computed by recursion and memoization. We show in Algorithm \ref{recursive alg} the pseudo-code for the dynamic programming approach. We use intermediate states \verb|S| to store the set $S_i$ of vertices $S_i$ added up to the step $i$. The function \verb|FacetComplexityRecursive(P, S)| computes recursively the minimum of $\sum_j a_j$ over all the partial sequences of facets that contain all the vertices not in $S$. In particular, \verb|FacetComplexityRecursive|$(P,\emptyset)=\rho'$. The same state \verb|S| may be encountered for many times along different paths. Hence, a memoization of those states can be faster. 

We implemented this algorithm as a {\em Macaulay2} function in the appended file \verb|code.m2|, which is also available through the author's homepage.

\begin{algorithm}[h]
\SetKwProg{Fn}{Function}{}{end}
\SetKwFunction{ComRecur}{FacetComplexityRecursive}
\caption{A recursive algorithm computing the facet complexity}\label{recursive alg}
\Fn{\ComRecur{P, S}}{
\KwData{A reflexive polytope $P$, an ordered list $S$ of vertices of $P$}
\KwResult{an integer $\rho'$}
initialization\;
$d \leftarrow \dim P$\;
\uIf{$\txt{length}(S) =\txt{length}(V(P))$ }{
\Return 0\;
}
\Else{
	$\rho' \leftarrow \infty$\;
	\For{$F\in F(P)$}{
		\If{ $V(F)\not\subseteq S$}{
		$r \leftarrow\txt{affineDim}(V(F)\cap S)$\;	
		$c \leftarrow \infty$\;
		\If{ $r \geq d-2$ }{
			$S' \leftarrow\txt{sort}(S\cup V(F))$\;
			$c \leftarrow d-1-r + \ComRecur(P, S')$\;
		}
		$\rho' \leftarrow \min(\rho', c)$\;
		}
	}
	\Return $\rho'$\;
}
}
\BlankLine
$\ComRecur(P, \emptyset)$\;
\end{algorithm}

\section{Picard ranks of a polar pair}\label{PolarPair}
In this section we consider the sum of the Picard ranks of a Gorenstein toric Fano variety and its polar variety. We show the sum is bounded from above by $\min \{\vv{F(P)}, \vv{V(P)}\} -d +1$. This can be seen as a generalization of \cite[2.16]{Eikelberg1993}.

\begin{theorem}\label{PQE}
If $d\geq 3$, and $X=X_P$ and $Y=Y_Q$ are two Gorenstein toric Fano varieties of dimension $d$, such that $P=Q^*$. Then
\[\rho_X + \rho_Y \leq \min \{\vv{F(P)}, \vv{V(P)}\} -d +1.\]
The equality holds if and only if $P$ or $Q$ is simplicial, in which case $\rho_Y=1$ or $\rho_X=1$ respectively.
\end{theorem}

\subsection{A tale of two adjacent sequences}
We introduced the facet complexity $\rho'$ in the previous section. The goal here is to construct simultaneously two adjacent sequences $T$ and $L$ over $V(P)$ and $F(P)$ respectively. Equivalently, by Definition \ref{adjSeq}, $T$ gives rise to a dual adjacent sequence $T^*$ over $F(Q)$. Then we show that $\rho_{L,P} +\rho_{T^*,Q}$ is bounded from above.

We first discussion some variations on the definition of the facet complexity. In Definition \ref{fc}, we refer to the condition such that $a_i=0$ as the `staying condition' $(S)$ and the one for $a_i=1$ the `moving-up condition' $(U)$. That is, $G_i$ satisfies the condition
\begin{itemize}
\item[($S$),] if $S_{i-1}\cap V(G_i)$ has affine dimension $d-1$; or
\item[($U$),] if $S_{i-1}\cap V(G_i)$ has affine dimension $d-2$.
\end{itemize}

We introduce the following variations:

\begin{enumerate}
\item\begin{itemize}
\item[($S_1$)] if $G_i\adj^{F} G_j$ for some $j< i$, and there exists another vertex $v$ of $G_i$ in $S_{i-1}\backslash F$. 
\item[($U_1$)] if $G_i\adj^{F} G_j$ for some $j< i$, and $S_{i-1}\cap G_i=V(F)$.

\end{itemize}
\item \begin{itemize}
\item[($S_2$)] if $G_i$ is adjacent to two different $G_j$, $j< i$. 
\item[($U_2$)] if $G_i$ is adjacent to exactly one $G_j$, $j< i$.
\end{itemize}
\end{enumerate}

It is easy to see that $(S_2)\implies (S_1) \implies (S)$, and if $L$ is an adjacent sequence, then $(S_1)$ and $(S)$, $(U_1)$ and $(U)$ are equivalent respectively for every $i\geq 1$. We can define in parallel $\rho'_1$ and $\rho'_2$, where we replace the conditions $(S)$ and $(U)$ in Definition \ref{fc} with the corresponding variations above. By Lemma \ref{adjSeqE} and the argument above, we have $\rho'\leq \rho'_1\leq \rho'_2<\infty$. We will use the equivalent conditions $(S_1, U_1)$ to prove the theorem.

The dual descriptions of the conditions $(S_1, U_1)$ on a vertex sequence of $Q$ are as follows. Consider an adjacent vertex sequence $T=(v_0,v_1,\cdots, v_t)$ on $P$. We say $v_i$ ($i\geq 1$) satisfies conditions ($S_1^*$) or ($U_1^*$) if:

\begin{itemize}
\item[($S_1^*$)] if $\overline{v_i v_j}$ is an edge of $P$ for some $j< i$, and $v_i\sim^G v_k$, for a different $k\neq j$, $k<i$ and a facet $G$ not containing $v_j$. 
\item[($U_1^*$)] if $\overline{v_i v_j}$ is an edge of $P$ for some $j\leq i$, and no such $v_k$ and $G$ above exist. 
\end{itemize}

\pfof{Theorem \ref{PQE}}. For any adjacent vertex sequence $T=(v_0,v_1,\cdots, v_{t-1})$ over $V(P)$, where $t=\vv{V(P)}$, we construct a facet sequence $L$ over $F(P)$ as follows.

Firstly, at step $0$, we add all the facets of in $\nbd(v_0)$ to $P$, in a suitable order that we designate later. That is, all the facets of $P$ that contains $v_0$. At step 1, we add those facets in $\nbd(v_1)$ that have not been added to $L$, if there are any, in an order such that they form an adjacent sequence. Inductively, at step $i$, we add those facets in $\nbd(v_i)$ that have not been added to $L$ up to now. Denote by $L_i$ the  set of facets added to $L$ at step $i$.

We show that for each $i\geq 0$, there is a way to arrange the facets in $L_i$ so that $L$ is an adjacent sequence. To this aim, we need only to arrange the facets in $L_i$ in a way such that they are all adjacent to a previous facet. 

For $i=0$, note that $\nbd(v_0)$ is dual to $V(F_0)$, where $F_0$ is the facet of $Q$ dual to $v_0$. Hence by Lemma \ref{adjSeqE}(1) applied to $F_0$, there is an adjacent sequence over $V(F_0)$. 

For $i\geq 1$, since $T$ is an adjacent sequence, $\overline{v_j v_i}$ is an edge for some $v_j$, $j< i$. There is a facet $G_{i,j}$ containing both $v_j$ and $v_i$. By the construction, $G_{i,j}$ is contained in $L$ at step at most $j$. Apply Lemma \ref{adjSeqE}(3) to $\nbd(v_i)$ and $L_i$, and we can add the facets in $L_i$ to $L$, such that each is adjacent to a previous facet in $L$. Thus by induction, $L$ is an adjacent sequence. Finally, since $T$ contains all vertices of $P$, $L$ must contain all the facets of $P$.

Define $\rho_{L,P}$ and $\rho_{T^*,Q}$ over $P$ and $Q$ respectively. We write $\rho_L$ and $\rho_{T^*}$ for short. By the proof of Theorem \ref{finitecomplexity}, both $\rho_L$ and $\rho_{T^*}$ are finite numbers. Let $k_i:=\vv{L_i}$. For each step $i\geq 1$, let $L_i=\{G_{i,1},\cdots, G_{i,k_i}\}$, in the order added to $L$. 

There are two cases. If $k_i=0$, then no facets are added to $L$ at this step. Then for any facet $H \in \nbd(v_i)$, and any vertex $v$ of $H$, $v\in\{v_0,\cdots, v_{i-1}\}$. Note that $v\in V(H)$ for some $H\in \nbd(v_i)$ if and only if $v\sim v_i$. Thus $\near(v_i)\subseteq \{v_0,\cdots, v_{i-1}\}$.
Since $d\geq 3$, there exist two vertices $w_1, w_2$ such that $v_i \adj w_1$ and $v_i \adj w_2$. Since $\adjt(v_i)\subseteq\near(v_i)$, we have $w_j=v_{i_s}$, for $i_s\leq i-1$, $s=1,2$. Now $v_i$ satisfies condition ($S_1^*$), because there must exists a facet $K$ containing the edge $\overline{v_i w_2}$ but not $w_1$. As a conclusion, the step $i$ contributes to the sum $\rho_L+\rho_{T^*}$ by zero.

Otherwise, if $k_i>0$, then each $G_{i,j}$ contributes to $\rho_L$ by either $0$ or $1$ for $1\leq j\leq k_i-1$. The last one $G_{i,k_i}$ contributes by $0$, because when $d\geq 3$, $G_{i,k_i}$ is adjacent to two other facets containing $v_i$, which have been added to $L$, and therefore $G_{i,k_i}$ satisfies $(S_2)$. The vertex $v_i$ contributes to $\rho_{T^*}$ by $0$ or $1$. As a result, the total contribution of $v_i$ and the facets in $L_i$ to the sum $\rho_L+\rho_{T^*}$ is at most $k_i$.

Consider the first vertex $v_0$ and all the $e=k_0$ neighboring facets $G_0,\cdots, G_{e-1}$ of $v_0$. Then $e=\vv{\nbd(v_0)}$ equals to the valence of $v_0$. By definition, $v_0$ doesn't contribute to $\rho_{T^*}$ as the initial vertex. Let $a_i$ be the contribution of $G_i$ to $\rho_L$ as defined in Definition \ref{fc}. By Theorem \ref{fc1}, we find that 
\begin{align*}
\rho_X + \rho_Y &\leq \rho_L + \rho_{T^*} \\
&\leq \sum_{j=1}^{e-1} a_j + \sum_{i=1}^r k_i\\
& =\sum_{j=1}^{e-1} a_j + \vv{F(P)}- e.
\end{align*}
By Lemma \ref{starA}, when $d\geq 3$, we can rearrange these $\{G_i\}_{i=0}^{e-1}$, so that $\sum_{j=1}^{e-1} a_j\leq e - (d-1)$. This shows 
\begin{equation}\label{FP}
\rho_X + \rho_Y\leq \vv{F(P)}- d + 1.
\end{equation}
By symmetry, we can replace $P$ with $Q$ and this proves that $\rho_X + \rho_Y\leq \min\{\vv{F(Q)},\vv{F(P)}\}- d + 1$. Then the inequality of the theorem follows since $\vv{F(Q)}=\vv{V(P)}$.

Finally we prove the condition when the equality holds. Suppose $Q$ is simplicial, then $\rho_Y = \vv{V(Q)}-d = \vv{F(P)}-d$. Besides, $\rho_X\geq 1$, as $X$ is projective. Hence $\rho_X +\rho_Y\geq \vv{F(P)}-d+1$. Therefore, $\rho_X +\rho_Y= \vv{F(P)}-d+1$. In particular, $\rho_Y= 1$.

Conversely, assume $\vv{F(P)}\leq \vv{V(P)}$ and the equality (\ref{FP}) holds. Then for any adjacent sequence $T$ of $V(P)$, there exists an adjacent sequence $L$ constructed as in the previous paragraph, such that $\rho_L + \rho_{T^*}= \vv{F(P)}- d + 1$. By the analysis above for the case $k_i>0$, we conclude that for every step $i\geq 1$, the contribution of $v_i$ in $T$ to $\rho_Y$ is one if and only if $k_i\geq 1$. 

Suppose $Q$ is not simplicial. Then there is a non-simplex facet $H$ of $Q$. Firstly $d\geq 3$. By Lemma \ref{miss2}, there is a facet $\sigma$ of $H$ which misses at least two vertices $x,y$ of $H$. Furthermore, there is a facet $\tau$ of $H$ containing $x$ but not $y$, since for any convex polytope, the intersection of all facets containing a vertex $v$ is $\{v\}$. In particular, there are two distinct facets $F_\sigma$, $F_\tau$ of $Q$ such that $F_\sigma \adj^{\sigma} H$, and $F_\tau \adj^{\tau} H$. Let $\mathcal{F}_y$ be the subset of $\adjt(H)$ that does not contain $y$. Then $F_\sigma, F_\tau\in \mathcal{F}_y$.

Over $Q$ we can construct an adjacent sequence $T^*=(F_0, \cdots, F_i=H, \cdots, F_m)$, such that $(F_0,\cdots, F_{i-1})$ is an adjacent sequence over $\mathcal{F}_y$. This follows from Lemma \ref{adjSeqE}(2). We have $F_i=H$ satisfies condition $(S_2)$, since $H$ is adjacent to two previous facets $F_\sigma$ and $F_\tau$. Thus $H$ contributes to $\rho_{T^*}$ by $0$. However, for any $L$ we constructed, $k_i\geq 1$, since the facet $F_y$ dual to $y$ is not added to $L$ until step $i$. So we reached a contradiction. This proves that $Q$ is simplicial.
\qed

\begin{lemma}\label{starA}
Suppose $P$ is a lattice polytope of dimension $d\geq 3$. Let $v\in V(P)$ and $F_0, \cdots, F_{e-1}$ be all the facets of $P$ containing $v$. Then there exists an adjacent sequence $L$ over $\{F_i\}$ such that $\sum_{j=1}^{e-1} a_j\leq e - (d-1)$.
\end{lemma}
\pf. Let $Q=P^*$, then the statement is equivalent to the following: for any facet $G$ of $Q$, there is an adjacent sequence $v_0,\cdots, v_{e-1}$ of $V(G)$, such that the number of $v_i$ satisfying condition ($U_1^*$) on $Q$ is at most $e - d + 1$.

Note any $v_i$ satisfies the condition ($S_1^*$) if and only if $v_i$ satisfies the following condition:

\begin{itemize}
\item[($X$)] if $\overline{v_i v_j}$ is an edge of $G$ for some $j< i$, and there exists a facet $\tau$ of $G$ and $k<i$, $k\neq j$, such that $v_j\not\in\tau$, and $v_i\sim^\tau v_k$ as faces of $G$. 
\end{itemize}
Thus we can restrict to the facet $G$. Then the lemma follows from Lemma \ref{starB}. \qed

\begin{lemma}\label{starB}
Suppose $G$ is a lattice polytope of dimension $d\geq 2$. Let $e=\vv{V(G)}$. Then there is an adjacent sequence $v_0,\cdots, v_{e-1}$ of $V(G)$, such that the number of $v_i$ ($i\geq 1$) not satisfying condition $(X)$ is at most $e-d$.
\end{lemma}
\pf. Firstly, $e\geq d+1$. For $d=2$ this is true for all $e\geq 3$, by choosing $v_0,\cdots, v_{e-1}$ in sequence along the boundary.

We prove the lemma by an induction on $k:=e-d$. When $k=1$, $G$ is a simplex, and the only vertex $v_i$ not satisfying condition ($X$) is $v_1$. Therefore the lemma holds for any adjacent sequence of vertices.

Suppose the lemma holds for all $k\leq k_0$. Consider the case when $G$ is of dimension $d$ with $e_0+1=d+k_0+1$ vertices. 
If $d=2$, then the lemma holds for $G$. Assume $d\geq 3$. By Lemma \ref{miss2}, there exists a facet $H$ of $G$ such that $H$ misses at least two vertices of $G$. Then $H$ is of dimension $d-1$ and has $e'\leq e_0 +1 -2 = d+k_0-1$ vertices. Now $H$ satisfies the assumption of the induction hypothesis, so there exists an adjacent sequence $L_H:=\{v_0,\cdots, v_{e'-1}\}$ over $V(H)$, such that the number of $v_i$ ($i\geq 1$) not satisfying condition (X) is at most $e'-d+1$. We now add to $L_H$ all the other vertices of $G$. By Lemma \ref{adjSeqE}(4), $L_H$ can be completed into an adjacent sequence $L_G$ over $V(G)$. Next, if $v_i$ satisfies condition (X) on $H$, then it satisfies condition (X) on $G$ too, since for any vertex $u$ of $H$, and any facet $\tau$ of $H$ such that $u\not\in\tau$, there exists a facet $F$ of $G$ such that $F\adj^\tau H$. Then $F$ containing $\tau$ but not $u$.

Since the last vertex in $L_G$ satisfies the condition (X), we find the the number of vertices in $L_G$, except for $v_0$, that do not satisfy condition (X), is at most $(e_0+1-e'-1)+(e'-d+1)=e_0+1-d$. By induction on $k$, this proves the lemma.\qed

\begin{lemma}\label{miss2}
Let $G$ be a convex polytope of dimension $d\geq 2$. If $G$ is not a simplex, then there exists a facet $H$ of $G$ such that $H$ misses at least two vertices of $G$.
\end{lemma}

\pf. Let $e:=\vv{V(G)}$. We show by induction that if a $d$-polytope $G$ satisfies that all the facets $H$ of $G$ misses exactly one vertex, then $G$ is a simplex, or equivalently, $e=d+1$.

For $d=2$ this holds, since any facet $H$ is a line segment, so if $H$ misses exactly one vertex, then $e=3$.
Suppose the claim holds for $d=d_0$. Assume $G$ is a $(d_0+1)$-polytope such that all its facets has exactly $e-1$ vertices. Any $(d_0-1)$-face $K$ of $G$ equals $H\cap H'$ for two facets $H, H'$ of $G$. Note $V(H\cup H') = V(G)$ for $d_0+1\geq 2$. Thus $\vv{V(K)}= \vv{V(H)}+\vv{V(H')}- \vv{V(H\cup H')}=2(e-1)-e=e-2$. Thus any facet $H$ of $G$ satisfies the inductive hypothesis, so $H$ is a $d_0$-simplex. Hence $e=d_0+2$, so $G$ is a simplex too. By induction this claim holds for all $d\geq 2$.\qed 

\begin{remark}\label{lbt}
When $P$ is simple, $Q$ is simplicial. Then $\rho_X=1$ and $\rho_Y=\vv{F(P)}-d$. Theorem \ref{PQE} implies that $\vv{V(Q)}\leq \vv{F(Q)}$. This in fact follows from the Lower Bound Theorem \cite{Barnette1973}: If $E$ is a simplicial $d$-polytope, then $\vv{F(E)}\geq (d-1)\vv{V(E)} - (d+1)(d-2)$. Since $\vv{V(E)}\geq d+1$, we have $\vv{F(E)}-\vv{V(E)}\geq (d-2)(\vv{V(E)}-(d+1))\geq 0$.
\end{remark}

\subsection{The difference $\varepsilon$}\label{Evalue}
We conclude by some investigations of the difference in Theorem \ref{PQE}. We define for every dual pair $(P,Q)$ of reflexive polytopes an integer
\[ \varepsilon = \varepsilon(P,Q):= \min \{\vv{F(P)}, \vv{V(P)}\} -d +1 - \rho_X - \rho_Y.\]
Then Theorem \ref{PQE} says $\varepsilon\geq 0$, with $\varepsilon=0$ if and only if $P$ or $Q$ is simplicial (so the other is simple).

\begin{example}\label{4309-2}
Consider the $3$-dimensional reflexive polytope in Example \ref{4309}. All the facets are parallelograms, hence non-simplicial. 

We know that $\rho_X=2, \rho_Y=4, \vv{F(P)}=12$, and $\vv{V(P)}=14$. Hence $\varepsilon = 4$.
\end{example}

\begin{example}\label{3dvare}
In dimension $d=3$, by computing the Picard ranks in {\em Macaulay2}, we find $0\leq \varepsilon\leq 9$. There is a unique reflexive $3$-polytope $P=\PALP(4317)=\GRD(1943)$, with $\varepsilon= 9$, and $P\cong P^*$ is self-dual. To be specific, $\rho_X=\rho_Y=1$, and $\vv{V(P)}=\vv{F(P)} =13$. Note that this is one of the 22 examples in Example \ref{neqpr} where $\rho'=2>\rho=1$.
\end{example}

For a reflexive polytope $Q$, let $\rho_{max}(Q)$ be the maximum of $\rho_L$ over all the sequences $L$ over $V(P)$ such that $\rho_L<\infty$.

\begin{lemma}\label{lowerboundv}
Let $(P,Q)$ be a polar pair of reflexive polytopes of dimension $d$. Suppose $\vv{F(P)}\leq \vv{V(P)}$. Then 
\[\varepsilon(P,Q)\geq (\rho'_P-\rho_X) + (\rho_{max}(Q)-\rho_Y).\]

\end{lemma}

\pf. The proof of Theorem \ref{PQE} shows that for every adjacent sequence $T$ over $V(P)$, there exists an adjacent sequence $L$ over $F(P)$, such that $\rho_{T^*,Q}+\rho_{L,P} \leq \vv{F(P)}-d+1$. By definition, $\rho_L\geq \rho'_P$. In particular, take $T$ such that $\rho_{T^*} = \rho_{max}(Q)$. Then 
\begin{align*}
\varepsilon(P,Q)&= \vv{F(P)}-d+1 - \rho_X- \rho_Y\\
&\geq\rho_{T^*} + \rho_L - \rho_X- \rho_Y\\
&\geq\rho_{max}(Q)+ \rho'_P - \rho_X- \rho_Y.
\end{align*}\qed

\bibliography{mybib.bib}
\bibliographystyle{alpha}
\end{document}